\documentclass{amsart}
\usepackage{epsfig}
\usepackage{amsfonts,amssymb,amscd,amsmath,euscript,enumerate,verbatim,calc}
\theoremstyle{plain}

\newtheorem{theorem}{Theorem}[section]
\newtheorem{lemma}[theorem]{Lemma}
\newtheorem{proposition}[theorem]{Proposition}

\newtheorem{corollary}[theorem]{Corollary}
\theoremstyle{definition}
\newtheorem{definition}[theorem]{Definition}

\newtheorem{remark}[theorem]{Remark}

\def\Fq{{\mathbb F}_q}
\def\FF{{\mathbb F}}

\newcommand{\rank}{\operatorname{rank}}
\newcommand{\CPP}{\operatorname{CPP}}
\newcommand{\HPP}{\operatorname{HP}}
\newcommand{\PnF}{\operatorname{P}_n(F)}
\newcommand{\CPPnF}{\CPP_n(F)}
\newcommand{\HPPnF}{\HPP_n(F)}
\newcommand{\GL}{\operatorname{GL}}
\newcommand{\T}{\operatorname{T}}
\newcommand{\TGL}{\operatorname{TGL}}
\newcommand{\Hn}{{\operatorname{H}}_n(F)}
\newcommand{\Hnq}{{\operatorname{H}}_n(\Fq)}
\newcommand{\HnrF}{{\operatorname{H}}_n(r,F)}

\newcommand{\HnkF}{{\operatorname{H}}_n^{(k)}(F)}
\newcommand{\HkkF}{{\operatorname{H}}_k^{(k)}(F)}
\newcommand{\HnkrF}{{\operatorname{H}}_n^{(k)}(r,F)}
\newcommand{\HnzeroF}{{\operatorname{H}}_n^{(0)}(F)}
\newcommand{\HnzerorF}{{\operatorname{H}}_n^{(0)}(r,F)}
\newcommand{\Hnrq}{{\operatorname{H}}_n(r,\Fq)}
\newcommand{\Hnrminusq}{{\operatorname{H}}_n(r-1,\Fq)}
\newcommand{\Hnkq}{{\operatorname{H}}_n^{(k)}(\Fq)}
\newcommand{\Hnnq}{{\operatorname{H}}_n^{(n)}(\Fq)}
\newcommand{\Hnkrq}{{\operatorname{H}}_n^{(k)}(r,\Fq)}
\newcommand{\Hnzeroq}{{\operatorname{H}}_n^{(0)}(\Fq)}
\newcommand{\Hnzerorq}{{\operatorname{H}}_n^{(0)}(r,\Fq)}
\newcommand{\HGL}{\operatorname{HGL}}

\newcommand{\CPPn}{\CPP_n(\Fq)}

\makeatletter
\def\Ddots{\mathinner{\mkern1mu\raise\p@
\vbox{\kern7\p@\hbox{.}}\mkern2mu
\raise4\p@\hbox{.}\mkern2mu\raise7\p@\hbox{.}\mkern1mu}}
\makeatother
\begin{document}
\title[Relatively Prime Polynomials and Nonsingular Hankel Matrices]{Relatively Prime Polynomials and Nonsingular Hankel Matrices over Finite Fields}

\author {Mario Garc\'{i}a-Armas}
\address{Department of Mathematics, University of British Columbia, \newline \indent 
Vancouver, BC V6T 1Z2, Canada.}
\email{marioga@math.ubc.ca} 

\author{Sudhir R. Ghorpade}
\address{Department of Mathematics, 
Indian Institute of Technology Bombay,\newline \indent
Powai, Mumbai 400076, India.}
\email{srg@math.iitb.ac.in}

\author{Samrith Ram}
\address{Department of Mathematics,
Indian Institute of Technology Bombay,\newline \indent
Powai, Mumbai 400076, India.}
\email{samrithram@gmail.com}



\begin{abstract}
The probability for two monic polynomials of a positive degree $n$ with coefficients in the finite field $\Fq$ 
to be relatively prime turns out to be identical with the probability for an $n \times n$ Hankel matrix over $\Fq$ to be nonsingular. 
Motivated by this, we give an explicit map from pairs of coprime polynomials to nonsingular Hankel matrices that explains this connection. 
A basic tool used here is the classical notion of Bezoutian of two polynomials. Moreover, we give simpler and direct proofs of the general 
formulae for the number of $m$-tuples of relatively prime polynomials  
over $\Fq$ of given degrees and for the number of $n\times n$ Hankel matrices over $\Fq$ of a given rank.
\end{abstract}
\date{\today}
\maketitle
\section{Introduction}
It is a remarkable fact that the probability for two randomly chosen monic polynomials of the same positive degree with coefficients in the binary field $\FF_2$ to be coprime is exactly $1/2$. This observation appears to go back at least to an exercise in the treatise, first published in 1969, 
of Knuth \cite[\S 4.6.1, Ex. 5]{K} (see also Remark~\ref{story}). More recently, it was made by Corteel, Savage, Wilf, and Zeilberger \cite{CSWZ} in 1998 in 
the course of their work on Euler's pentagonal sieve in the theory of partitions, and it led them to ask for a ``nice simple bijection'' between the 
coprime and the non-coprime ordered pairs of monic polynomials of degree $n$ over $\FF_2$. This was answered first by Reifegerste \cite{Reif} in 2000 and by Benjamin and Bennett \cite{BB} in 2007. The latter deals with the more general case of polynomials over any finite field $\Fq$ where the probability turns out to be $1-(1/q)$ instead of $1/2$. Since there are $q^{2n}$ ordered pairs of monic polynomials over $\Fq$ of degree $n$, this means~that
\begin{equation}
\label{cpp} 
\left|\CPPn\right| = q^{2n} - q^{2n-1} = q^{2n-1}(q-1),
\end{equation}
where $\CPPn$ denotes the set of ordered pairs of coprime monic polynomials over $\Fq$ of degree $n$. In effect, Benjamin and Bennett gave an explicit surjective map from $\CPPn$ onto the set of ordered pairs of non-coprime monic polynomials over $\Fq$ of degree $n$ in such a way that the cardinality of 
each fiber is $q-1$. 

A couple of years prior to \cite{CSWZ} and working on a seemingly unrelated topic, 
Kaltofen and Lobo \cite{KL} observed that the probability for a $n\times n$ Toeplitz matrix with entries in $\Fq$ to be nonsingular is exactly $1-(1/q)$. In fact, this observation can be traced back to Daykin \cite{D} who had essentially proved the same result (and also a more general one) in 1960 with Hankel matrices
in place of Toeplitz matrices. 
Since there are $q^{2n-1}$  Toeplitz matrices (or equivalently, Hankel matrices) of size $n\times n$ with entries in $\Fq$, this means that
\begin{equation}
\label{tgl} 
\left|\TGL_n(\Fq)\right| = \left|\HGL_n(\Fq)\right| = q^{2n-1} - q^{2n-2} = q^{2n-2}(q-1),
\end{equation}
where $\TGL_n(\Fq)$ (resp: $\HGL_n(\Fq)$) denotes the set of all $n\times n$ nonsingular Toeplitz (resp: Hankel) matrices with entries in $\Fq$. 

One of the main aims of this paper is to explain the uncanny coincidence 
that the probability in both of the above situations turns out to be the same or, more precisely, the fact that 
the formulae \eqref{cpp} and \eqref{tgl} differ just by a factor of $q$. We do this by giving an explicit surjective map from 
$\CPPn$ onto $\HGL_n(\Fq)$ such that each fiber has cardinality $q$. This readily yields a similar map with $\HGL_n(\Fq)$ replaced by $\TGL_n(\Fq)$. As a consequence, we obtain new proofs of \eqref{cpp} and \eqref{tgl} by combining any one of the known proofs with this surjective map. We further add to this collection of proofs by giving  alternative, short and completely self-contained  proofs of more general versions of \eqref{cpp} and \eqref{tgl}. 
%

\section{Preliminaries} 
\label{prelim}


Let $F$ be a field. 
Recall that a matrix $M=\left(m_{ij}\right)$ with entries in $F$ is said to be a \emph{Toeplitz matrix} (resp: \emph{Hankel matrix}) if $m_{ij} = m_{rs}$ 
whenever $i-j=r-s$ (resp: $i+j=r+s$). Thus every $n\times n$ Toeplitz (resp: Hankel) matrix   
over $F$ looks like $\left(a_{n+i-j}\right)$ (resp: $\left(a_{i+j-1}\right)$)  
for a unique $(2n-1)$-tuple $\left(a_1,  \dots, a_{2n-1}\right)\in F^{2n-1}$.

We denote by $\T_n(F)$ (resp: $\Hn$) the set of all Toeplitz (resp: Hankel) matrices with entries in $F$ and, as in the Introduction, 
set 
$$
\TGL_n(F) = \T_n(F) \cap \GL_n(F) \quad \text{ and } \quad  
\HGL_n(F) = \Hn \cap \GL_n(F).
$$

The following simple observation shows that at least as far as enumerative and bijective combinatorics is concerned, Toeplitz and Hankel matrices are the same.  

\begin{proposition}
\label{TandH}
There is a bijection between $\T_n(F)$ and $\Hn$, which induces a bijection between $\TGL_n(F)$ and $\HGL_n(F)$.
\end{proposition}

\begin{proof}
Let $E$ be 
the $n\times n$ matrix with $1$ on the antidiagonal and $0$ elsewhere, i.e., $E=\left(\delta_{i,n-j+1}\right)$ 
where $\delta$ is the Kronecker delta. Then $E$ is nonsingular and the map given by $A\mapsto AE$ sets up the desired bijection.
\end{proof}

As usual, $F[X]$ will denote the set of polynomials in one variable $X$ with coefficients in $F$.  
Recall that for any $u,v\in F[X]$ 
of degree $\le n$, 
the $n$th order \emph{Bezoutian} (matrix) of $u$ and $v$ is 
the $n\times n$ matrix $B_n(u,v)= \left(b_{ij}\right)$ determined by the 
equation
$$
\frac{u(X)v(Y) - v(X)u(Y)}{X-Y} = \sum_{i,j=1}^n b_{ij} X^{i-1}Y^{j-1}.
$$
The coefficients $b_{ij}$ are not hard to determine explicitly; in fact, if 
$u = \sum_{i=0}^n u_i X^i$ and $v = \sum_{i=0}^n v_i X^i$, then 
upon letting $u_k=v_k:=0$ for $k>n$, we have 
$$
b_{ij} = \sum_{s=1}^{\min\{i, j\}} \left(v_{s-1} u_{i+j-s} - u_{s-1} v_{i+j-s}\right) \quad \text{ for $1\le i,j\le n$} .
$$
It is clear from the definition that if $u$ and $v$ have a nonconstant common factor then the system of homogeneous linear equations corresponding to $B_n(u,v)$ has a nontrivial solution, and hence $B_n(u,v)$ is singular. It is a classical fact that the converse is also true; we record this below for convenience and refer to the survey article 
\cite{HF} of Helmke and  Fuhrmann for a proof. 

\begin{proposition}
\label{bez}
Let $u, v\in F[X]$. 
Assume that $\deg u =n$ and $\deg v \le n$. 
Then  $B_n(u,v)$ is nonsingular if and only if 
$u$ and $v$ are coprime. 
\end{proposition}

As an illustration, consider $u,v\in F[X]$ such that $v$ is the constant polynomial~$1$ and 
$u(X) = u_0+u_1X+\dots + u_nX^n$ with $u_0, u_1, \dots , u_n\in F$.  
Then 
$$
\frac{u(X) - u(Y)}{X-Y} = 
\sum_{k=1}^n u_k\frac{X^k - Y^k}{X-Y} = \sum_{k=1}^n u_k\sum_{i=1}^k X^{i-1}\, Y^{k-i}
= \sum_{i,j=1}^n u_{i+j-1} X^{i-1}\, Y^{j-1},
$$
where, by convention, $u_k:=0$ for $k>n$. Thus the $n$th order Bezoutian $B_n(u,1)$ has $u_n$ on its antidiagonal and $0$ below that. In particular, if $\deg u=n$, i.e., if $u_n\ne 0$, then 
$u$ and $v$ are coprime, and moreover $B_n(u,v)$ is nonsingular.  

\section{An Explicit Surjection} 
\label{surjmap}

Fix a positive integer $n$ and a field $F$. 
As in the Introduction, let 
$$
\CPPnF : = \left\{(f,g)\in F[X]^2 : f, \text{$g$ are coprime and both are monic of degree $n$}\right\}.
$$
Moreover, let us consider
\begin{align*}
\PnF & : = \left\{(u,v)\in F[X]^2 :  \text{$u$ is monic, $\deg u =n$, and } \deg v < n  \right\}, \text{ and}  \\
\HPPnF  & : = \left\{(u,v)\in \PnF : u \text{ and $v$ are coprime} \right\}.
\end{align*}
 
We may refer to an element of $\PnF$ as a \emph{Pad\'e pair} and an element of $\HPPnF$ 
as a \emph{Hermite pair}. 
\begin{lemma}
\label{fgtouv}
$\CPPnF$ is in bijection with $\HPPnF$.
\end{lemma}

\begin{proof} 
The map given by $(f,g)\mapsto (f, g-f)$ does the job. 
\end{proof}

\begin{lemma}
\label{pade}
Let $(u,v) \in \PnF$. Then there are unique $a_i\in F$, $i\ge 1$, such that 
\begin{equation}
\label{uva}
\frac{v(X)}{u(X)} = \sum_{i=1}^{\infty} \frac{a_i}{X^i}.
\end{equation}
\end{lemma}

\begin{proof}
Write $u(X) = X^n\left[1 - u^*(1/X)\right]$ for a unique 
$u^*\in F[X]$ with no constant term. 
Expanding as a formal power series, we obtain 
$$
\frac{v(X)}{u(X)} = 
X^{-n} v(X)\sum_{j=0}^{\infty} u^*(1/X)^j.
$$
This yields the desired $a_i\in F$.   
\end{proof}

\begin{definition}
For $(u,v) \in \PnF$, we define $H_n(u,v)$ to be the $n\times n$ Hankel 
matrix whose $(i,j)$th entry is $a_{i+j-1}$ for $1\le i,j\le n$, 
where $a_1, a_2, \dots$ are as in \eqref{uva}.
\end{definition}

The following result which relates $H_n(u,v)$ to the $n$th order Bezoutian $B_n(u,v)$ 
is classical and is sometimes referred to as Barnett's factorization. 
We include a proof for the sake of completeness, especially since the proofs found in the literature are often a bit involved and 
tend have an additional assumption that the polynomials 
$u$ and $v$ are coprime, i.e., $(u,v)$ is a Hermite pair rather than a Pad\'e pair. 

\begin{proposition}
\label{barnett}
$B_n(u,v) = B_n(u,1)H_n(u,v)B_n(u,1)$ for any $(u,v) \in \PnF$. 
\end{proposition}

\begin{proof}
Let $R(T):=v(T)/u(T)$ and let $a_i$, $i\ge 1$, be as in \eqref{uva}. Then 
$$
\frac{R(Y)-R(X)}{X-Y} = \sum_{i=1}^{\infty} a_i \sum_{j=1}^i \frac{X^{i-j}\, Y^{j-1}}{X^iY^i} = \sum_{k,\ell=1}^{\infty} a_{k+\ell-1}X^{-k}\, Y^{-\ell}.
$$
Now if $u(X) = u_0+\dots + u_{n-1}X^{n-1} + X^n$ with $u_0,  \dots , u_{n-1}\in F$ and  $u_n:=1$, then 
\begin{align*}
\frac{u(X)v(Y) - v(X)u(Y)}{X-Y} &= u(X) \frac{R(Y)-R(X)}{X-Y} u(Y) \\
&= \Bigg(\sum_{r=0}^nu_rX^r\Bigg)\Bigg(\sum_{k,\ell=1}^{\infty} a_{k+\ell-1}X^{-k}Y^{-\ell}\Bigg)\Bigg(\sum_{s=0}^n u_sY^s\Bigg) \\ 
&=  \sum_{i,j\le n} \Bigg(\sum_{k,\ell \ge 1}u_{i+k-1}\, a_{k+\ell-1}\, u_{\ell+j-1}\Bigg)X^{i-1}Y^{j-1},
\end{align*}
where, by convention, $u_t=0$ for $t>n$ and $a_t=0$ for $t\le 0$. Comparing the coefficients of $X^{i-1}Y^{j-1}$ for $1\le i,j\le n$, we obtain the desired result. 
\end{proof}

\begin{theorem}
There is a surjective map $\sigma:\CPPnF \to \TGL_n(F)$ such that for any $A\in \TGL_n(F)$, the fiber $\sigma^{-1}\left(\{A\}\right)$ is in one-to-one correspondence with $F$. In particular, 
$
\left|\CPPn \right|= q \left|\TGL_n(\Fq)\right|.
$
\end{theorem}
 
\begin{proof}
{F}rom Propositions \ref{bez} and \ref{barnett}, we see that 
$H_n(u,v)$ is nonsingular for any $(u,v)\in \HPPnF$. Thus, we obtain a well-defined map $\eta:\HPPnF \to \HGL_n(F)$ given by $(u,v)\mapsto  H_n(u,v)$. 
Now let $B\in \HGL_n(F)$. Then there are unique $b_1, \dots , b_{2n-1}\in F$ such that the $(i,j)$th entry of $B$ is $b_{i+j-1}$ for $1	\le i,j \le n$.
Let $\lambda$ be an arbitrary element of $F$ and set $b_{2n}:=\lambda$. Since $B$ is nonsingular, there are unique $u_0,  \dots, u_{n-1}\in F$ such that 
\begin{equation}
\label{Bua}
B
\begin{pmatrix} 
u_0  \\ u_{1} \\ \vdots \\ u_{n-1}
\end{pmatrix}
= 
- \begin{pmatrix} 
b_{n+1}  \\ b_{n+2} \\ \vdots \\ b_{2n}
\end{pmatrix}.
\end{equation}
Next, define $u_n:=1$ and $v_0, v_1, \dots , v_{n-1}$ to be the unique elements of $F$ given by the following triangular system of equations:
\begin{equation}
\label{Buv}
\begin{pmatrix} 
v_0  \\ v_{1} \\ \vdots \\ v_{n-1}
\end{pmatrix}
= \begin{pmatrix} 
b_1  & b_{2} & \dots  & b_n \\
0  & b_{1} & \dots  & b_{n-1} \\
\vdots &  &\ddots &  \vdots \\
0  & 0 & \dots  & b_1 \\
\end{pmatrix}
 \begin{pmatrix} 
u_{1}  \\ u_{2} \\ \vdots \\ u_{n}
\end{pmatrix}.
\end{equation}
Finally, define $u, v\in F[X]$ by 
$$
u =  \sum_{i=0}^n u_i X^i   \quad \text{ and } \quad v = \sum_{i=0}^{n-1} v_i X^i.
$$
Then $(u,v)\in \PnF$ and if we let $a_k\in F$, $k\ge 1$, as in \eqref{uva}, then we have 
$$
\sum_{i=1}^{n} v_{i-1} X^{i-1} = \Bigg(\sum_{j=0}^n u_j X^j \Bigg)\Bigg(\sum_{k\ge 1} a_k X^{-k} \Bigg) = \sum_{i\le n} \Bigg(\sum_{j=0}^n a_{j-i+1}u_j \Bigg)X^{i-1},
$$
where, by convention, $a_k:=0$ for $k\le 0$. Comparing the coefficients of $X^{i-1}$ for $-n< i\le n$, we find that \eqref{Buv} and \eqref{Bua}
are satisfied with
$b_1, \dots , b_{2n}$ replaced by $a_1, \dots , a_{2n}$, respectively. Since $u_n=1$, the triangular nature of \eqref{Buv} implies that 
$a_i = b_i$ for $1\le i\le n$. Further, successive comparison of \eqref{Bua} with its counterpart where $b_i$'s are replaced by $a_i$'s yields 
$a_i = b_i$ for $1\le i\le 2n$.  In particular, $B=H_n(u,v)$. Now since $B$ is nonsingular, Propositions \ref{barnett} and \ref{bez} show that $u$ and $v$ are coprime. Thus $(u,v)\in \HPPnF$ and $\eta(u,v) = B$. It is clear from the construction above 
that a Hermite pair $(u,v)$ satisfying $\eta(u,v) = B$ 
is uniquely determined by the matrix 
$B$ and the element $b_{2n}=\lambda$. 
Also, in view of \eqref{Bua}, distinct values of $\lambda$ 
in $F$ give rise to distinct monic polynomials $u$ in $F[X]$ of degree $n$. This shows that  for each $B\in \HGL_n(F)$, the fiber $\eta^{-1}\left(\{B\}\right)$ is in one-to-one correspondence with $F$. Finally, combining $\eta$ with the bijections given by Proposition \ref{TandH} and Lemma \ref{fgtouv}, we obtain the desired surjective map $\sigma:\CPPnF \to \TGL_n(F)$. 
\end{proof}
  
\section{Relatively Prime Polynomials}

The general version of \eqref{cpp} alluded to in the Introduction is the theorem stated below. 
It may be noted that this generalizes 
\cite[Prop. 3]{CSWZ}, \cite[Thm. 9]{N} and \cite[Prop. 2.4]{MZ}, and 
also that it is a more precise form of \cite[Cor. 5]{BB} 
and \cite[Thm. 1.1]{HM}. 
We remark at the outset that in this theorem, considering arbitrary polynomials (not necessarily monic) in $\Fq[X]$ 
does not affect the probability. 

\begin{theorem}
\label{probmn}
Let $m$ be a positive integer and $n_1, \dots , n_m$ be nonnegative integers. The probability that $m$ monic 
polynomials in $\Fq[X]$ of degrees $n_1, \dots , n_m$, chosen independently and uniformly at random, are relatively prime is $1 - q^{1-m}$ if $\min\{n_1, \dots , n_m\}\ge 1$ and $1$ otherwise. 
\end{theorem} 

\begin{proof} 
Let $N(n_1, \dots , n_m)$ denote the number of ordered $m$-tuples 
$\left(f_1, \ldots, f_m\right)$ of coprime monic polynomials in $\Fq[X]$ such that $\deg f_i = n_i$ for $i = 1, \ldots, m$.
Evidently, it suffices to show that 
\begin{equation}
{N(n_1, \ldots,n_m)} = \begin{cases}
q^{n_1+ \cdots+n_m} \left(1-q^{1-m}\right) & \textrm{ if $\min\{n_1, n_2, \ldots,n_m\} \geq 1$},\\
q^{n_1+ \cdots+n_m}                       & \textrm{ if $\min\{n_1, n_2, \ldots,n_m\} = 0$}.
\end{cases}
\end{equation}
To this end, we shall assume, without loss of generality, that $n_1 \geq  \dots \geq n_m$. 
We can partition the set 
of ordered $m$-tuples $(f_1, \dots , f_m)$ of monic polynomials in $\Fq[X]$ with $\deg f_i =n_i$ for $i\le i\le m$, into disjoint subsets
$S_0, S_1, \ldots, S_{n_m}$, where   for $0\le d\le n_m$, 
the set $S_d$ consists of $m$-tuples whose GCD is of degree $d$.
Given any monic polynomial $h\in \Fq[X]$ of degree $d$ and any coprime $m$-tuple $\bigl(g_1, \ldots, g_m\bigr)$ of monic polynomials such that $\deg g_i = n_i - d$ for $i = 1, \ldots, m$, it is easy to see that $\left(hg_1, \dots , hg_m\right)\in S_d$. 
Conversely, if $\bigl(f_1, \ldots, f_m\bigr) \in S_d$, then the polynomial
$h = {\rm GCD}\bigl(f_1, \ldots, f_m\bigr)$ is monic of degree $d$ and $(f_1/h, \dots , f_m/h)$
is an ordered $m$-tuple of coprime monic polynomials of degrees $n_1-d, \ldots, n_m-d$, respectively.  
This shows that $\left|S_d\right| = q^{d} N(n_1-d, \ldots, n_m-d)$ for $0\le d\le n_m$, and consequently, 
\begin{equation} \label{Eq1}
q^{n_1+ \cdots+n_m} = \sum_{d=0}^{n_m} \left|S_d \right|= \sum_{d=0}^{n_m} q^{d} N(n_1-d, n_2-d, \ldots, n_m-d).
\end{equation}
If $n_m = 0$, we immediately obtain $N(n_1, \ldots, n_m) = q^{n_1+\cdots+n_m}$. On the other hand, if $n_m \geq 1$, substituting $n_i$ by $n_i-1$ ($i=1, \ldots, m$) in the above relation yields
\begin{equation}\label{Eq2}
q^{n_1+ \cdots+n_m - m} = \sum_{d=1}^{n_m} q^{d-1} N(n_1-d, n_2-d, \ldots, n_m-d).
\end{equation}
Multiplying equation \eqref{Eq2} by $q$ and subtracting the result from \eqref{Eq1}, we obtain 
$N(n_1,n_2, \ldots,n_m) = q^{n_1+\cdots+n_m} \left(1 - q^{1-m}\right)$, as desired. 
\end{proof}

\begin{remark} 
\label{story}
As indicated in the Introduction, the case $m=2$ (and $q$ prime) of the above result appears as an exercise (\# 5 of \S 4.6.1) in Knuth \cite{K}.
The solution outlined by Knuth uses the result obtained in the previous exercise and in turn, a deep analysis of the euclidean algorithm. 
The general result with arbitrary $m$ and $q$ (but with $n_1=\dots = n_m=n$) given in Corteel, Savage, Wilf and Zeilberger \cite{CSWZ} seems to 
have been arrived at independently by completely different means. Also, it is indicated in a footnote in \cite[p.188]{CSWZ} that the degrees $n_1, \dots , n_m$ could 
well be different (i.e., one has a result such as Theorem \ref{probmn} above), and this observation is ascribed to D. Zagier. Many of the subsequent works (e.g. \cite{BB,HM,GP,Reif}) cite \cite{CSWZ} as an earliest reference for this result 
(and in fact, the authors of this paper did the same before it was pointed out by
a referee that the result is classical). In retrospect, the key ideas in the answer by Benjamin and Bennett \cite{BB} to the question 
in \cite{CSWZ} about a nice bijective proof can be traced back to \cite[\S 4.6]{K} and a more detailed analysis by Norton \cite{N} as well as by 
Ma and van zur Gathen \cite{MZ}. In the same vein, the short proof given above of Theorem \ref{probmn}, even though it was discovered independently,  
can be viewed as an extension of the ``alternative proof'' that appears in the solution of Exercise 5 of \S 4.6.1 in the first edition of Knuth \cite{K}, but for some mysterious reason, is missing in the subsequent editions. Thus, the contents of this section may help resurrect an original  and perhaps the simplest proof. Finally, we remark that nontrivial generalizations of Theorem \ref{probmn} are studied by Gao and Panario \cite{GP} and by Hou and Mullen \cite{HM}, while an application to a conjecture about the enumeration of certain Singer cycles is discussed in \cite{GS}.
\end{remark}

\section{Hankel Matrices over $\Fq$}

The general version of \eqref{tgl} alluded to in the Introduction is the following.

\begin{theorem}
\label{Hnr}
The number $N(n,r; q)$ of $\, n\times n$ Hankel matrices of rank $\, r$ with entries in the finite field $\Fq$ is given by 
\begin{equation}
\label{Nnrq}
N(n,r; q) = \begin{cases}
1 & \text{ if } \; r=0,\\
q^{2r-2}(q^2-1) & \text{ if } \; 1 \leq r \leq n-1,\\
q^{2n-2}(q-1) & \text{ if } \; r=n.
\end{cases}
\end{equation}
\end{theorem} 

Before giving a proof of the above theorem, we introduce some notation and prove a few auxiliary results. 
Let $F$ be a field and, as before, $n$ a positive integer. Given any $n\times n$ matrix $A$ with entries in $F$ and any positive integers $d$,  
$i_1, \dots , i_d, j_1, \dots , j_d$ such that $i_1< \dots < i_d \le n$ and $j_1< \dots < j_d \le n$, we denote by $A[i_1, \dots , i_d | j_1, \dots , j_d]$ the $d\times d$ submatrix of $A$ formed by the rows indexed by $i_1, \dots , i_d$ and the columns indexed by $j_1, \dots , j_d$. 
Note that the $d^{\rm th}$ leading principal submatrix of $A$ is $A[1, \dots , d|1, \dots , d]$ and this will be denoted simply by $A_d$. Define
\begin{equation*}
\delta(A) := \begin{cases}
0 & \text{ if $A_d$ is singular for each }  d=1,\dots , n,\\
\max\{d: \text{$A_d$ is nonsingular}\} & \text{ otherwise.}
\end{cases} 
\end{equation*}
For 
$r,k\in \{0,1, \dots ,n\}$, let $\HnrF :=\left\{A\in \Hn : \rank(A) \le r \right\}$, and moreover, 
$$
\HnkF : = \left\{A\in \Hn : \delta(A) = k \right\} \quad\text{and} \quad
\HnkrF:= \HnrF\cap \HnkF. 
$$
Note that $\HGL_n(F)={{\operatorname{H}}_n^{(n)}(F)}={{\operatorname{H}}_n^{(n)}(n,F)}$ and also that 
\begin{equation}
\label{HnDecomp}
\Hn = \coprod_{k=0}^n \HnkF \quad\text{and} \quad \HnrF = \coprod_{k=0}^r \HnkrF,
\end{equation}
where $\coprod$ denotes disjoint union. The main idea in the proof of Theorem \ref{Hnr} 
is to use the above decompositions and 
to characterize $\Hnkq$ and $\Hnkrq$ suitably so as to be able to determine their cardinalities recursively. 
Here is the first step. 

\begin{lemma}
\label{lem2_1}
Let $A=\left(a_{i+j-1}\right)\in \Hn$. 
Then 
\begin{equation}
\label{Hnzero}
A\in \HnzeroF \Longleftrightarrow a_1=\dots = a_n =0.
\end{equation}  
Moreover, for $0\le r\le n-1$,
\begin{equation}
\label{Hnzeror} 
A\in \HnzerorF \Longleftrightarrow a_1=\dots = a_{2n-r-1} =0.
\end{equation}
In particular, $|\Hnzeroq| = q^{n-1}$ and $|\Hnzerorq| = q^{r}$.
\end{lemma}

\begin{proof}
If $A\in \HnzeroF$, then $\det(A_k)=0$ for $k=1, \dots ,n$. Using this successively, we obtain $a_1=\dots = a_n =0$. Conversely, if $a_1=\dots = a_n =0$, then it is clear that  $\det(A_k)=0$ for $k=1, \dots ,n$, i.e., $A\in \HnzeroF$. Next, let $0\le r\le n-1$  
and suppose $A\in \HnzerorF$. Then 
$a_1=\dots = a_n =0$, as before. Moreover, by successively using the vanishing of the $(r+1)\times (r+1)$ minor $\det A[n-r, n-r+1, \dots , n|j, j+1, \dots , j+r]$ for $j=2, \dots , n-r$, we obtain $a_{n+1}= \dots = a_{2n-r-1}=0$ as well. Conversely, suppose 
$a_1=\dots = a_{2n-r-1} =0$, then  $A\in \HnzeroF$ and it is easily seen that every $(r+1)\times (r+1)$ submatrix of $A$ has a column of zeros, and so $A\in \HnzerorF$. 
\end{proof}

The following result is an analogue of \eqref{Hnzero} for $\HnkF$ where $k\ge 1$. 

\begin{lemma}
\label{lem2_2}
Let $k\in \{1, \dots ,n-1\}$ and $A=\left(a_{i+j-1}\right)\in \Hn$ be such that $A_k$ is nonsingular. Suppose $\mathbf{x}=\left(x_1, \dots , x_k\right)\in F^k$ is the unique solution of the system 
$A_k\mathbf{x}^T = \left(a_{k+1}, \dots , a_{2k}\right)^T$, i.e., 
for $t=1, \dots , k$, the following relation holds: 
\begin{equation}
\label{lincombn} 
a_{k+t} = x_1 a_{t} + \dots + x_k a_{t+k-1}. 
\end{equation}
Then 
\begin{equation}
\label{HnkIff}
A\in \HnkF \Longleftrightarrow \text{ the relation \eqref{lincombn} holds for $t=1, \dots , n$.}
\end{equation}
\end{lemma}

\begin{proof}
Suppose $A\in \HnkF$. 
We will use induction on $t$ to show that \eqref{lincombn} holds for $t=1,\dots , n$. The case when $1\le t\le k$ is known by the hypothesis. So let us assume that $t\ge k+1$ and that \eqref{lincombn} holds for all values of $t$ smaller than given one. Consider the $t\times t$ matrix $A_t$ and successively make the following $t-k$ elementary column transformations:
$$
C_t - \left(x_1 C_{t-k}  + \dots + x_kC_{t-1}\right), \; 
\dots , \; C_{k+1} - \left(x_1C_1+\dots +x_kC_k\right) 
$$
where $C_j$ indicates the $j^{\rm th}$ column. This transforms $A_t$ to the $t\times t$ matrix 
\begin{equation*}
A'=
\left( \begin{array} {ccc|ccc} 
a_1&\cdots&a_{k}&h_{1}&\cdots &h_{t-k}\\
\vdots&\Ddots&\vdots&\vdots& &\vdots\\
a_{k}&\cdots&a_{2k-1}&h_{k}&\cdots & h_{t-1}\\
\hline
a_{k+1}&\cdots&a_{2k}&h_{k+1}& \cdots&h_{t}\\
\vdots& &\vdots&\vdots&\Ddots&\vdots\\
a_{t}&\cdots&a_{k+t-1}&h_{t}&\cdots &h_{2t-k-1}
\end{array} \right),
\end{equation*}
where $h_m = a_{k+m} - (x_1a_m + \dots + x_k a_{m+k-1})$ for $m=1, \dots 2t-k-1$. By induction hypothesis, $h_m=0$ for $m=1,\dots , t-1$, and therefore
$$
\det(A_t) = \det(A') = (-1)^{(t-k)(t-k+1)/2} \det(A_k) h_t^{t-k}.
$$
Since $\det(A_t)=0$ and $\det(A_k)\ne 0$, it follows that $h_t=0$, i.e., the relation \eqref{lincombn} holds for the given value of $t$. 

Conversely, suppose the relation \eqref{lincombn} holds for $t=1, \dots , n$. Then we can write 
${\mathbf{v}}_{k+1}= x_1{\mathbf{v}}_1+ \dots + x_k{\mathbf{v}}_k$, where ${\mathbf{v}}_j$ denotes the $j^{\rm th}$ column vector of $A$. In particular, $\rank(A_{k+1})\le k$, which implies that $A\in \HnkF$.
\end{proof}

Let us pause to observe that the formula \eqref{tgl} for the number of nonsingular Hankel matrices can already be derived as a 
consequence of the above results. 

\begin{corollary}
\label{NonsingHankel}
$|\Hnzeroq| = q^{n-1}$ and $\big|\Hnkq\big| = q^{n+k-2}(q-1)$ for $1\le k \le n$. In particular, $\left|\HGL_n(\Fq)\right|= q^{2n-2}(q-1)$.
\end{corollary}

\begin{proof}
Induct on $n$. If $n=1$, then $k$ is $0$ or $1$, and the desired formulae are obvious. Suppose $n>1$ and the result holds for positive values of $n$ smaller than the given one. By Lemma \ref{lem2_1},  $|\Hnzeroq| = q^{n-1}$. Now suppose $1\le k < n$.  Then by Lemma~\ref{lem2_2}, we see that the map 
$A=\left(a_{i+j-1}\right) \mapsto \left(A_k, \; a_{2k}, a_{n+k+1}, \dots , a_{2n-1}\right)$ gives a bijection of $\HnkF$ onto $\HkkF \times F^{n-k}$. 
Hence using the induction hypothesis,  $\big|\Hnkq\big| = q^{2k-2}(q-1)q^{n-k} = q^{n+k-2}(q-1)$. Finally, in view of \eqref{HnDecomp},  
the induction hypothesis, and an easy evaluation of a telescopic sum, we conclude that 
$
\left|\HGL_n(\Fq)\right|= 
\big|\Hnnq\big| = \big|\Hnq\big| - \sum_{k=0}^{n-1} \big|\Hnkq\big| = q^{2n-1} - q^{2n-2}.
$
\end{proof} 

If a Hankel matrix in $\HnkF$ satisfies a rank condition, the validity of \eqref{lincombn}  can be pushed a little further. More precisely, one has the following analogue of \eqref{Hnzeror}. 

\begin{lemma}
\label{lem2_3}
Let $k,r$ be integers with 
$1\le k \le r < n$ 
and $A=\left(a_{i+j-1}\right)\in \Hn$ be such that $A_k$ is nonsingular.  
Suppose $\mathbf{x}=\left(x_1, \dots , x_k\right)\in F^k$ is the unique solution of the system 
$A_k\mathbf{x}^T = \left(a_{k+1}, \dots , a_{2k}\right)^T$.  
Then 
\begin{equation}
\label{HnkrIff}
A \in \HnkrF \Longleftrightarrow \text{ the relation \eqref{lincombn} holds for $t=1, \dots , 2n-r-1$.}
\end{equation} 
\end{lemma}

\begin{proof}
Suppose $A \in \HnkrF$. 
Again, we use induction on $t$. By Lemma \ref{lem2_2}, the relation \eqref{lincombn} holds if $1\le t \le n$. Assume that 
$n+1\le t\le 2n-r-1$ 
and that \eqref{lincombn} holds for all values of $t$ smaller than the given one. Define $\mathbf{x}^{(0)}, \mathbf{x}^{(1)}, \dots , \mathbf{x}^{(t-1)}$ in $F^k$ recursively as follows. First, $\mathbf{x}^{(0)} :=\mathbf{x} =   \left(x_1, \dots , x_k\right)$. Next, if ${\ell}\ge 1$ and if 
$\mathbf{x}^{({\ell}-1)}= \big(x^{({\ell}-1)}_1, \dots , x^{({\ell}-1)}_k\big)$ is known, then we let $x_0^{({\ell}-1)}:=0$ and let 
$\mathbf{x}^{({\ell})}= \big(x^{({\ell})}_1, \dots , x^{({\ell})}_k\big)\in F^k$ 
be given by 
$$
x^{({\ell})}_s = x_s x^{({\ell}-1)}_k + x^{({\ell}-1)}_{s-1} \quad \text{ for } s=1, \dots , k. 
$$
Observe that for $1\le {\ell} < t $ and $1\le m < t$, we have
\begin{equation}
\label{xjs}
\sum_{s=1}^k x^{({\ell})}_s a_{m+s-1} = x^{({\ell}-1)}_k\sum_{s=1}^k x_s a_{m+s-1} + \sum_{s=1}^{k-1}  x^{({\ell}-1)}_s a_{m+s} = \sum_{s=1}^{k}  x^{({\ell}-1)}_s a_{m+s}, 
\end{equation}
where the last equality follows from \eqref{lincombn} with $t$ replaced by $m$. Successive application of \eqref{xjs} shows that
\begin{equation}
\label{xjsxs}
\sum_{s=1}^k x^{({\ell})}_s a_{m+s-1} = \sum_{s=1}^k x_s a_{m+s+{\ell}-1} \quad \text{for } 0\le \ell < t 
\text{ and } 1\le m \le t-\ell. 
\end{equation}
Now consider the 
$(2n-t)\times (2n-t)$ principal submatrix $B$ of $A$ given by 
$$
B:= A[1, 2, \dots , k, t+k-n+1,  \dots, n-1,  n \vert 1, 2, \dots , k, t+k-n+1,  \dots, n-1, n]
$$ 
and make the following $2n-t-k$ elementary column transformations:
$$
C_{2n-t} - \sum_{s=1}^k x^{(n-k-1)}_s C_{s} , \; 
C_{2n-t-1} - \sum_{s=1}^k x^{(n-k-2)}_s C_{s} ,
\dots , \; C_{k+1} - \sum_{s=1}^k x^{(t-n)}_s C_{s} , 
$$
where $C_j$ indicates the $j^{\rm th}$ column. This transforms $B$ to the $(2n-t)\times (2n-t)$~matrix 
\begin{equation*}
B'=
\left( \begin{array} {ccc|ccc} 
a_1&\cdots&a_{k}&u_{1,\, 1}&\cdots &u_{1, \, 2n-t-k}\\
\vdots&\Ddots&\vdots&\vdots& &\vdots\\
a_{k}&\cdots&a_{2k-1}&u_{k,\, 1}&\cdots & u_{k, \, 2n-t-k}\\
\hline
a_{t+k-n+1}&\cdots&a_{t+2k-1}&v_{1,\, 1}& \cdots&v_{1,\, 2n-t-k}\\
\vdots& &\vdots&\vdots&\Ddots&\vdots\\
a_{n}&\cdots&a_{n+k-1}&v_{2n-t-k,\, 1}&\cdots &v_{2n-t-k, \, 2n-t-k}
\end{array} \right),
\end{equation*}
for some $u_{i,\, j}, \; v_{i,\, j}\in F$. In fact, for $1\le i\le k$ and $1\le j\le 2n-t-k$,
\begin{align*}
u_{i,\, j} & =  a_{i+t+k-n+j-1} - \sum_{s=1}^k x^{(t-n+j-1)}_s a_{i+s-1} \\
       & =   a_{i+j+t-n-1+k} - \sum_{s=1}^k x_s a_{i+j+t-n-1+s-1} = 0,
\end{align*}
where the penultimate equality follows from \eqref{xjsxs} and the last equality follows from \eqref{lincombn} since $t\ge n+1$. 
Moreover, for $1\le i, j\le 2n-t-k$,
$$
v_{i,\, j}  =  a_{2t+2k-2n+i+j-1} - \sum_{s=1}^k x^{(t-n+j-1)}_s a_{t+k-n+i+s-1};
$$
also, since $t\ge n+1$, using \eqref{xjsxs}, we have 
\begin{equation}
\label{vij}
v_{i,\, j}     =   a_{2t+2k-2n+i+j-1} - \! \sum_{s=1}^k x_s a_{2t+k-2n+i+j+s-2} \quad \text{if $i+j\le 2n-t-k+1$.}  
\end{equation}
In particular, $v_{i,\, j}$ depends only on $i+j$ whenever $i+j\le 2n-t-k+1$. 
Furthermore, if $i+j\le 2n-t-k$, then 
from \eqref{lincombn}  we deduce that $v_{i,\, j}=0$ .
Consequently, upon letting $v= v_{1,\,  2n-t-k} = \cdots  = v_{2n-t-k,\, 1}$, we obtain 
$$
\det(B) = \det(B') = (-1)^{(2n-t-k)(2n-t-k+1)/2}  \det(A_k) v^{2n-t-k}. 
$$
But since $A\in \HnkrF$ and $2n-t \ge r+1$, we have 
$\det(A_k)\ne 0$ and $\det(B)=0$. Hence $v=0$,  and from \eqref{vij}, 
we conclude that \eqref{lincombn} holds for the given value of $t$. 

Conversely, suppose the relation \eqref{lincombn} holds for $t=1, \dots , 2n-r-1$. 
Then we can write 
${\mathbf{v}}_{j}= x_1{\mathbf{v}}_{j-k}+ \dots + x_k{\mathbf{v}}_{j-1}$ for $j=k+1, \dots, k+n-r$,  where ${\mathbf{v}}_j$ denotes the $j^{\rm th}$ column vector of $A$. Hence the column space of $A$ is 
spanned by ${\mathbf{v}}_1, \dots , {\mathbf{v}}_k, {\mathbf{v}}_{k+n-r+1}, \dots , {\mathbf{v}}_n$. 
In particular, 
$\rank(A) \le k+ n-(k+n-r+1)+1 = r$.
This together with Lemma \ref{lem2_2} shows that $A\in \HnkrF$.
\end{proof}

\begin{corollary}
\label{hglnandhnkrq}
$|\Hnzerorq| = q^{r}$  for $0\le r < n$ and $\big|\Hnkrq\big| = q^{r+k-1}(q-1)$ for $1 \leq k \le r < n$. 
Consequently, $\left|\Hnrq\right| = q^{2r}$ for $0\le r < n$. 
\end{corollary} 

\begin{proof}
The first assertion follows from Lemma \ref{lem2_1}. Now suppose $1 \leq k \le r < n$. 
By Lemma~\ref{lem2_3}, we see that the map 
$A=\left(a_{i+j-1}\right) \mapsto \left(A_k, \; a_{2k}, a_{2n-r+k}, \dots , a_{2n-1}\right)$ gives a bijection of $\HnkrF$ onto $\HGL_k(F) \times F^{r-k+1}$. 
Hence by Corollary \ref{NonsingHankel},  $\big|\Hnkq\big| = q^{2k-2}(q-1)q^{r-k+1} = q^{r+k-1}(q-1)$. Finally, $\left|\Hnrq\right| = q^{2r}$ is obvious when $r=0$, whereas if $1\le r<n$, then using \eqref{HnDecomp}  
 and an easy evaluation of a telescopic sum, we conclude that 
$\big|\Hnrq\big| =  \sum_{k=0}^{r} \big|\Hnkrq\big| = q^{2r}$.
\end{proof} 

We are now ready to prove the main result of this section. 

\medskip

\noindent
\emph{Proof of Theorem \ref{Hnr}:} The case $r=0$ is trivial 
and for $r=n$, Corollary \ref{NonsingHankel} applies. %
Finally, if $1\le r<n$, then 
 $N(n,r;q)= \left|\Hnrq\right| - \left|\Hnrminusq\right| = q^{2r} - q^{2r-2}$,  
thanks to Corollary \ref{hglnandhnkrq}. 
\hfill $\Box$

\medskip

A noteworthy consequence of Theorem \ref{Hnr} is that for a fixed positive integer $r$, the number of $n\times n$ Hankel matrices of rank $r$ remains 
constant for every $n\ge r+1$. 

\section*{Acknowledgments}
  
We are grateful to Bharath Sethuraman for his help in bringing together the two sets of authors from two different continents. We also thank an anonymous referee for bringing 
\cite[\S 4.6]{K} and \cite{GP} to our attention.


\end{document}